**Probabilistic analytical approach to determining the asymptotics of prime objects on the initial interval of the natural series**

Victor Volfson

ABSTRACT. This paper considers a probabilistic-analytical approach to determining asymptotics of prime objects on the initial interval of the natural series. The author proposes a new method based on the construction of a probability space.

An arithmetic function is considered that counts the number of prime objects (for example, prime numbers, twin primes, prime values of polynomials) on a given probability space. The properties of this arithmetic function are proved. Based on these properties, the specified probabilistic-analytical approach is constructed.

As examples of the application of this approach, the definition of the asymptotics of the number of prime twins and pairs of prime numbers that add up to an even number (based on Goldbach's conjecture) is considered. It is shown that the proposed method allows one to obtain asymptotic estimates that coincide with the known conjecturea of prime number theory.

This approach opens up new possibilities for studying conjectures about prime numbers, offering an alternative way to prove them based on probabilistic methods.

# 1. INTRODUCTION

Prime number theory is one of the most problematic areas in number theory. The number of conjectures about prime numbers is steadily increasing, while proofs are lacking. Here are only the main ones in historical order: Goldbach's conjecture on the representation of even numbers greater than 2 as the sum of two primes [1], Riemann's conjecture on the accuracy of the number of primes on the initial interval of the natural series [2], Hardy-Littlewood conjecture on the number of prime tuples on the initial interval of the natural series [3], Cramer's conjecture on the maximum distance between consecutive primes on the initial interval of the natural series [4], Bateman-Horne conjecture on the number of prime values of polynomials on the initial interval of the natural series [5]. By the initial interval of the natural series, we mean the segment $[2, N]$, where $N$ is a sufficiently large natural number.

Among the main scientific results obtained, it is important to note the proof of the prime number theorem (PNT) on the number of prime numbers [6] and for the arithmetic progression [7] with accuracy that is not as high compared to the Riemann conjecture. The other conjectures mentioned above remain unproven. What is the matter?

The proof of PNT was obtained using analytical methods of number theory, which is based on the multiplicative structure of primes. However, when studying the asymptotics of primes, this approach does not justify itself. The point is that although primes have a lot of explicit multiplicative structure, on the other hand, the product of two primes is never prime. The revealed statistical anomalies in the asymptotics of primes cannot be easily explained in terms of the multiplicative properties of primes. Thus, when considering atypical statistics about primes, they behave pseudo-randomly and, therefore, can be approximated with acceptable accuracy using probabilistic models.

The most accurate are random probability models for prime numbers that include random variables. This is despite the fact that prime numbers are clearly deterministic in nature. Even for multiplicative problems, which are in principle controlled by the zeros of the Riemann zeta function, one can get good predictions by assuming various pseudorandom properties of these zeros, so that the distribution of these zeros can be modeled using a random model.

Of course, one cannot expect absolute accuracy when copying a deterministic set such as the prime numbers using a probabilistic model of that set, and each of the heuristic models discussed below has limitations on the range of statistics about the prime numbers that they can be expected to track with reasonable accuracy.

## 1.1 Cramer's random model

Cramer in 1936 [4] presented a probabilistic model of prime numbers with the following assumptions:

1. Prime numbers are distributed "randomly". The probability that a randomly chosen number $n$ is prime is approximately $1/\ln(n)$.

2. Prime numbers are "independent" of each other in the sense, for example, that the probability of two numbers both $n$ and $n+2$ being prime is equal to the product of the probabilities of each being prime.

3. To estimate the number of combinations of prime numbers up to some number $N$, we sum the probabilities for all possible combinations of prime numbers, such as combinations $(n, n+2)$ for twin primes.

This random model allowed Cramer to make an asymptotic estimate for the maximum distance between consecutive primes with probability 1 (almost everywhere). Cramer's random model allows one to prove the Riemann conjecture, as well as the Legendre conjecture and others, almost everywhere.

Despite its usefulness, Cramer's model has several well-known shortcomings. The probability of an even natural number $n$ being prime is not $1/\ln(n)$, but is 0. Similarly for natural numbers of the form $n = mq$, where $q$ is a prime and $m$ is a natural number. In addition, the assumption that primes occur independently of each other, such as prime twins, is not true.

## 1.2 Improved Cramer's model

To correct the above-mentioned shortcomings of Cramer's random model, Grenville [4] modified this model as follows. He successively discarded natural numbers multiples of prime numbers: $2, 3, ..., q$ on the initial interval of the natural series. Only natural numbers coprime with $Q = 2 \cdot 3 \cdot ... \cdot q$ remained after that. Thus, natural numbers of the form $n = mq$ were discarded, which gave an error in the probability estimate.

This model captures the correct global distribution of primes. However, unlike Cramer's model, the improved model also captures the prime bias in the residue classes modulo primes. In particular, this model satisfies the Hardy-Littlewood conjecture for counting the number of prime tuples [3].

The improved Cramer model allowed Grenville to refine the conjecture about the asymptotic estimate of the maximum distance between consecutive primes [8].

A common drawback of probabilistic models of prime numbers is that they allow one to prove conjectures about prime numbers based on some unproven assumptions and only with probability equal to 1 (almost everywhere).

We will not make assumptions of the indicated probability models in this work and will rely only on proven facts.

The aim of this work is attempt to develop a new probabilistic analytical approach to proving asymptotic conjectures about prime numbers.

Let us consider the following probability space.

Any initial segment of the natural series $\{1,2,...,n\}$ can be naturally transformed into a discrete probability space $\left(\Omega_n, \mathcal{A}_n, \mathbb{P}_n\right)$ by taking $\Omega_n = \{1,2,...,n\}$, $\mathcal{A}_n$ — all subsets $\Omega_n$, $P_n(A) = \#(m \in A)/n$, where $\#(m \in A)$ is the number of natural numbers in the subset [9].

Then an arbitrary arithmetic function $f(m), m = 1,...,n$ on $\Omega_n$ can be considered as a random variable $x_n$ on this probability space:

$$x_n(m) = f(m)(1 \le m \le \mathrm{n}).$$

The number of objects of prime numbers not exceeding a natural number $n$ is a real arithmetic function, so it can be considered as a random variable on the specified probability space. Let us denote this arithmetic function as $K(n)$.

The number of of prime tuples and the number of prime values of polynomials have already been considered in the implementation of this probabilistic analytical approach to the proof of the Hardy-Littlewood and Bateman-Horn conjectures in [10].

We will look at the properties of arithmetic function $K(n)$ in the next chapter.

## 2. PROPERTIES OF ARITHMETIC FUNCTION $K(n)$.

We will consider the arithmetic function $K(n)$, which can be represented as $K(n) = \sum_{i=1}^{n} 1_A(i)$, where $1_A(i)$ is the indicator function.

Assertion 1

The asymptotic behavior of the sum $\sum_{i=1}^{n} K(i)/i$ will not change when discarding a finite number of first terms.

Proof

$K(n)=\sum_{i=1}^{n} 1_A(i)$ is a non-decreasing arithmetic function with values in the interval: $1 \le K(n) \le n$, therefore $1/n \le K(n)/n \le 1$. It follows that $\sum_{n=1}^{\infty} K(n)/n \ge \sum_{n=1}^{\infty} 1/n$. Since the series $\sum_{n=1}^{\infty} 1/n$ diverges, the series $\sum_{n=1}^{\infty} K(n)/n$ diverges.

Let us denote the sum of a finite number of first terms of the series - $\sum_{i=1}^{b} K(i)/i = B$, then the limit of the ratio:

$$\lim_{n\to\infty} \frac{\sum_{i=b+1}^{n} K(i)/i}{\sum_{i=1}^{n} K(i)/i} = \lim_{n\to\infty} \frac{\sum_{i=b+1}^{n} K(i)/i}{B+\sum_{i=b+1}^{n} K(i)/i} = \lim_{n\to\infty} \frac{1}{1+B/\sum_{i=b+1}^{n} K(i)/i} = 1,$$

since the series $\sum_{i=b+1}^{\infty} K(i)/i$ - diverges. The assertion is proved.

Assertion 2

Let there is the arithmetic function $K(n)$ correspond to an equally probable random variable with values $K(1), K(2), ..., K(n)$ on the above-mentioned probability space constructed on the initial interval of the natural series $(\Omega_n, \mathcal{A}_n, \mathbb{P}_n)$. Let the probabilities on the spaces from $(\Omega_1, \mathcal{A}_1, \mathbb{P}_1)$ to $(\Omega_n, \mathcal{A}_n, \mathbb{P}_n)$ also be given: $p_1 = K(1)/1, ..., p_n = K(n)/n$.

Then the asymptotic is true $K(n) \sim \sum_{i=3}^{n} K(i)/i = \sum_{i=3}^{n} p_i$ for $K(i) = Ci/\ln^k(i)$ and $n \to \infty$, where $C$ is a real number and $k$ is a natural number.

Proof

Having in mind that the relation $\sum_{i=2}^{n} C/\ln^k(i) \sim Cn/\ln^k(n)$ is true for $n \to \infty$, then the asymptotics $K(n) \sim \sum_{i=3}^{n} K(i)/i = \sum_{i=3}^{n} p_i$ is true for $K(i) = Ci/\ln^k(i)$ and $n \to \infty$, where $C$ is a real number and $k$ is a natural number.

Corollary 3

Based on assertions 1 and 2 the asymptotic is true for $K(i) = Ci/\ln^k(i)$ ($C$ is a real number and $k$ is a natural number) and $n \to \infty$:

$$K(n) \sim \sum_{i=3}^{n} K(i)/i \sim \sum_{i=b+1}^{n} K(i)/i = \sum_{i=b+1}^{n} p_i .$$

These properties of the arithmetic function $K(n)$ are the basis of the probabilistic analytical approach, which will be discussed in the next chapter.

## 3. PROBALISTIC ANALYTICAL APPROACH TO DETERMINING THE ASYMPTOTICS OF PRIME OBJECTS ON THE INITIAL INTERVAL OF NATURAL SERIES

The probability of a large natural number $n$ to be prime, based on PNT, is approximately:

$$P(n) \approx 1/\ln n . \tag{3.1}$$

This is not quite true. Suppose that $n \geq 4$ is even, then the specified probability should be $0$. Similarly for cases when $n$ is a multiple of prime numbers: $2,3,...,q$.

Let us eliminate the indicated drawback. Let there are natural numbers for which (3.1) is true. Let us remove from them the numbers that are divisible by $2,3,...,q$. Now we denote $Q = \prod_{2 \leq p \leq q} p$, and the subset of the natural series obtained after removing the numbers - $A$ .

Let $a \in A$, then based on construction - $\gcd(a,Q) = 1$. Thus, in total there are suitable $\varphi(Q)$ residue classes.

The class of residues modulo $Q$ forms an arithmetic progression;

$$a, a+Q, a+2Q, ... \quad \text{c } \gcd(a,Q) = 1 \tag{3.2}$$

Let $n \in A$, using the theorem on the distribution of prime numbers in arithmetic progressions, based on (3.2), the number of prime numbers in this arithmetic progression when $n \to \infty$ is asymptotically equal to:

$$\pi(Q,n) \sim \frac{n}{\varphi(Q)\ln n}.$$

Since the number of natural numbers in the arithmetic progression (3.2) is asymptotically equal to $n/Q$, the asymptotics of the density of prime numbers in the subset $A$ obtained after removing composite numbers multiples of $2,3,...,q$ increases and, for $n \to \infty$, is equal to:

$$d(Q,n) = \frac{\pi(Q,n)}{n/Q} \sim \frac{Q}{\varphi(Q)\ln n}. \qquad (3.3)$$

Based on (3.3), the asymptotics of the probability (see note below) for $n \to \infty$ is equal to:

$$P(n) \sim \frac{Q}{\varphi(Q)\ln n}. \qquad (3.4)$$

Note. Let $B_n$ is a finite subset of the natural numbers with a sufficiently regular structure (allowing us to apply methods of number theory and probability theory) and with a number of terms tending to infinity as $n \to \infty$. Then we can construct a probability space similar to that constructed earlier in Chapter 1 of the work. $\Omega_n = B_n$, $\mathcal{A}_n$ are all subsets of $\Omega_n$, $P_n(A) = \#(m \in A)/\#(m \in B_n)$ in this case, i.e. equal to the density $A$ on $B_n$, and the arithmetic function $f(m)$ on $\Omega_n$ can be considered as a random variable $x_n(m) = f(m),(m \in B_n)$. The arithmetic progression (3.2) is used as $B_n$ in formula (3.4).

Let there is a sequence of natural numbers $a_1,...,a_k$ satisfying the conditions $\gcd(a_1,Q)=1,...,\gcd(a_k,Q)=1$. Let also $n >> Q$ and $a_1 > n,...,a_k > n$ then the different events that $a_i \in A$, will be asymptotically independent, i.e. when $n \to \infty$ the condition is satisfied:

$$P(a_1 \in A,...,a_k \in A) = P(a_1 \in A)...P(a_k \in A) . \qquad (3.5)$$

The reason for the asymptotic independence of these events is that $a_i \in A$ is not divisible by any small prime number up to $q$ (according to Buchstab's small prime number theorem [11]).

Thus, taking into account (3.4) and (3.5) we obtain:

$$P(a_1 \in A, ..., a_k \in A) = P(a_1 \in A)...P(a_k \in A) \sim \frac{Q^k}{\varphi(Q)^k \ln^k(n)}. \quad (3.6)$$

Based on (3.6) and Corollary 3:

$$K(N) \sim \sum_{n=Q}^{N} \frac{Q^k}{\varphi(Q)^k \ln^k n}. \quad (3.7)$$

where the lower bound is justified by the fact that the condition $\gcd(a,Q)=1$ may be violated at $n<Q$.

Based on (3.7), the asymptotic is equal to at $n \to \infty$:

$$K(N) \sim \frac{Q^k}{\varphi(Q)^k} \int_{t=Q}^{N} \frac{dt}{\ln^k t} = C \int_{t=Q}^{N} \frac{dt}{\ln^k t}. \quad (3.8)$$

The values $C,k$ are the real and natural numbers in (3.8) depend on the specific conjecture about prime numbers.

4. EXAMPLES OF THE PROBALISTIC ANALITICAL APPROACH

As an example of the indicated probabilistic analytical approach, we consider the definition of the asymptotic of the number of twin primes on the initial interval of the natural series $[2,N]$, where $N$ is a large number.

Based on (3.7), the following asymptotics is true:

$$K(N) \sim \sum_{n=Q}^{N} \frac{Q^2}{\varphi(Q)^2} \frac{1}{\ln^2(n)}. \quad (4.1)$$

Based on the Chinese Remainder Theorem [12], there are $\prod_{3 \le p \le q} (p-2)$ classes of residues modulo $Q = \prod_{2 \le p \le q} p$, such that $\gcd(n,Q)=1, \gcd(n+2,Q)=1$, so taking into account (4.1) we obtain:

$$\sum_{n=Q}^{N} \frac{Q^2}{\varphi(Q)^2 \ln^2(n)} = [N/Q] \prod_{2 \le p \le q} (p-2) \frac{Q^2}{\varphi(Q)^2 \ln^2(n)}. \quad (4.2)$$

Since $\varphi(Q)=\varphi(\prod_{2\le p\le q} p)=\prod_{2\le p\le q}(p-1)$, having in mind (4.1) and (4.2), we obtain the desired asymptotics $K(N)$ at $N\to\infty, q\to\infty$:

$$K(N)\sim[N/Q]\prod_{2\le p\le q}(p-2)\frac{Q^2}{\varphi(Q)^2\ln^2(n)}\sim\frac{2N}{\ln^2(N)}\prod_{p\ge3}\frac{p(p-2)}{(p-1)^2}. \tag{4.3}$$

Let us consider another example in which we will determine the asymptotic of the number of pairs of prime numbers that add up to an even number that does not exceed a larger natural number $N$ (based on Goldbach's conjecture).

Let $a\in A$, $N-a\in A$, $a\le N/2$, where $\gcd(a,Q)=1$.

Having in mind the asymptotic independence of events $a_i\in A$, based on (3.6):

$$P(a\in A, N-a\in A)=P(a\in A)P(N-a\in A)\sim\frac{Q^2}{\varphi(Q)^2\ln^2(N)}. \tag{4.4}$$

Based on (3.7), the following asymptotics is true:

$$K(N/2)\sim\sum_{i=Q}^{N/2}\frac{Q^2}{\varphi(Q)^2\ln^2(i)}. \tag{4.5}$$

Based on the Chinese Remainder Theorem, there are $\prod_{p\backslash N}(p-2)\prod_{p/N}(p-1)$ classes of residues modulo $Q=\prod_{2\le p\le q} p$, such that $(a,Q)=1,(N-a,Q)=1$, where $p\backslash N$ - means that the product is taken over all primes that do not divide $N$. Therefore, taking into account (4.5), we obtain:

$$\sum_{i=Q}^{N/2}\frac{Q^2}{\varphi(Q)^2\ln^2(i)}=[N/2Q]\prod_{3\le p\le q}(p-2)\prod_{p/N}\frac{(p-1)}{(p-2)}\frac{Q^2}{\varphi(Q)^2\ln^2(N)}. \tag{4.6}$$

Since $\varphi(Q)=\varphi(\prod_{2\le p\le q} p)=\prod_{2\le p\le q}(p-1)$, then having in mind (4.6), we obtain the desired asymptotics $K(N/2)$ $N\to\infty, q\to\infty$:

$$K(N/2)\sim[N/2Q]\prod_{3\le p\le q}(p-2)\prod_{p/N}\frac{(p-1)}{(p-2)}\frac{Q^2}{\varphi(Q)^2\ln^2(N)}\sim\frac{N}{\ln^2(N)}\prod_{p\ge3,p/N}\frac{(p-1)}{(p-2)}\prod_{p\ge3}(1-\frac{1}{(p-1)^2}). \tag{4.7}$$

## 5. CONCLUSION AND SUGGESTIONS FOR FURTHER WORK

The article explores a new probabilistic-analytical approach to determining the asymptotics of prime objects on the initial interval of the natural series. In the future, it is interesting to use this approach to other conjectures about prime numbers.

## 6. ACKNOWLEDGEMENTS

Thanks to everyone who has contributed to the discussion of this paper. The author is very grateful to Professor Gérald Tenenbaum at the University of Lorraine, who showed interest in the article and read the full text before publishing it in the archive.